\documentclass[11pt,a4paper]{amsart}
\usepackage[colorlinks=true,citecolor=black,linkcolor=black]{hyperref}
\usepackage{amsthm,amsfonts,amsmath, amssymb,graphicx,float, color}
\usepackage[numbers,sort&compress]{natbib}
\usepackage[lmargin=33mm,rmargin=33mm,tmargin=33mm,bmargin=33mm]{geometry}
\newcommand{\doi}[1]{\href{http://dx.doi.org/#1}{doi:\,\texttt{#1}}}
\renewcommand{\baselinestretch}{1.21}
\theoremstyle{plain}
\newtheorem{theorem}{Theorem}[section]
\newtheorem{lemma}[theorem]{Lemma}
\newtheorem{corollary}[theorem]{Corollary}
\newtheorem{proposition}[theorem]{Proposition}
\newtheorem{claim}[theorem]{Claim}

\theoremstyle{definition}
\newtheorem{conjecture}[theorem]{Conjecture}

\DeclareMathOperator{\dist}{dist}

\newcommand{\eps}{\varepsilon}  
\newcommand{\ceil}[1]{\ensuremath{\protect\lceil{#1}\rceil}}
\newcommand{\floor}[1]{\ensuremath{\protect\lfloor{#1}\rfloor}}
\newcommand{\Ceil}[1]{\ensuremath{\protect\left\lceil{#1}\right\rceil}}

\newcommand{\CeilFrac}[2]{\Ceil{\frac{#1}{#2}}}
\newcommand{\ceilFrac}[2]{\ceil{\frac{#1}{#2}}}
\newcommand{\Brac}[1]{\ensuremath{\protect\left(#1\right)}}

\newcommand{\half}{\ensuremath{\tfrac12}}

\newcommand{\NumEdges}[1]{\ensuremath{\|{#1}\|}}
\newcommand{\md}[1]{$K_{#1}$-model}

\begin{document}

\title{Small Minors in Dense Graphs}

\thanks{Samuel Fiorini and Gwena\"el Joret are supported in part by
  the Actions de Recherche Concert\'ees (ARC) fund of the Communaut\'e
  fran\c{c}aise de Belgique. Gwena\"el Joret is a Postdoctoral
  Researcher of the Fonds National de la Recherche Scientifique
  (F.R.S.--FNRS). David Wood is supported by a QEII Research
  Fellowship from the Australian Research Council.}

\author{Samuel Fiorini} \address{\newline D\'epartement de
  Math\'ematique \newline Universit\'e Libre de Bruxelles \newline
  Brussels, Belgium} \email{sfiorini@ulb.ac.be}

\author{Gwena\"el Joret} \address{\newline D\'epartement
  d'Informatique \newline Universit\'e Libre de Bruxelles \newline
  Brussels, Belgium} \email{gjoret@ulb.ac.be}

\author{Dirk Oliver Theis} \address{\newline Fakult\"at f\"ur
  Mathematik \newline Otto-von-Guericke-Universit\"at Magdeburg
  \newline Magdeburg, Germany} \email{dirk.theis@ovgu.de}

\author{David R. Wood} \address{\newline Department of Mathematics and
  Statistics \newline The University of Melbourne \newline Melbourne,
  Australia} \email{woodd@unimelb.edu.au}


\date{May 6, 2010. Revised: \today}

\begin{abstract} A fundamental result in structural graph theory
  states that every graph with large average degree contains a large
  complete graph as a minor. We prove this result with the extra
  property that the minor is small with respect to the order of the
  whole graph. More precisely, we describe functions $f$ and $h$ such
  that every graph with $n$ vertices and average degree at least
  $f(t)$ contains a $K_t$-model with at most $h(t)\cdot\log n$
  vertices. The logarithmic dependence on $n$ is best possible (for
  fixed $t$). In general, we prove that $f(t)\leq 2^{t-1}+\eps$. For
  $t\leq 4$, we determine the least value of $f(t)$; in particular
  $f(3)=2+\eps$ and $f(4)=4+\eps$.  For $t\leq4$, we establish similar
  results for graphs embedded on surfaces, where the size of the
  $K_t$-model is bounded (for fixed $t$).
\end{abstract}

\maketitle

\section{Introduction}

A fundamental result in structural graph theory states that every
sufficiently dense graph contains a large complete graph as a
minor\footnote{We consider simple, finite, undirected graphs $G$ with
  vertex set $V(G)$ and edge set $E(G)$. Let $|G|:=|V(G)|$ and
  $\NumEdges{G}:=|E(G)|$. A graph $H$ is a \emph{minor} of a graph $G$
  if $H$ is isomorphic to a graph obtained from a subgraph of $G$ by
  contracting edges.}. More precisely, there is a minimum function
$f(t)$ such that every graph with average degree at least $f(t)$
contains a $K_t$-minor. \citet{Mader67} first proved that $f(t)\leq
2^{t-2}$, and later proved that $f(t)\in O(t\log t)$
\citep{Mader68}. \citet{Kostochka82,Kostochka84} and
\citet{Thomason84,Thomason01} proved that $f(t)\in\Theta(t\sqrt{\log
  t})$; see \citep{Thomason06} for a survey of related results.


Here we prove similar results with the extra property that the
$K_t$-minor is `small' with respect to the order of the graph.  This
idea is evident when $t=3$. A graph contains a $K_3$-minor if and only
if it contains a cycle. Every graph with average degree at least $2$
contains a cycle, whereas every graph $G$ with average degree at least
$3$ contains a cycle of length $O(\log|G|)$. That is, high average
degree forces a short cycle, which can be thought of as a small
$K_3$-minor.

In general, we measure the size of a $K_t$-minor via the following
definition. A \emph{$K_t$-model} in a graph $G$ consists of $t$
connected subgraphs $B_1,\dots,B_t$ of $G$, such that $V(B_i)\cap
V(B_j)=\varnothing$ and some vertex in $B_i$ is adjacent to some
vertex in $B_j$ for all distinct $i,j\in\{1,\dots,t\}$. The $B_i$ are
called \emph{branch sets}. Clearly a graph contains a $K_t$-minor if
and only if it contains a $K_t$-model. We measure the size of a
$K_t$-model by the total number of vertices, $\sum_{i=1}^t|B_i|$. Our
main result states that every sufficiently dense graph contains a
small model of a complete graph.

\begin{theorem}
  \label{th:Main}
  There are functions $f$ and $h$ such that every graph $G$ with
  average degree at least $f(t)$ contains a $K_t$-model with at most
  $h(t)\cdot\log|G|$ vertices.
\end{theorem}

For fixed $t$, the logarithmic upper bound in Theorem~\ref{th:Main} is
within a constant factor of being optimal, since every $K_t$-model
contains a cycle, and for all $d \geq 3$ and $n>3d$ such that $nd$ is
even, \citet{Chandran03} constructed a graph with $n$ vertices,
average degree $d$, and girth at least $(\log_d n)-1$.  (The
\emph{girth} of a graph is the length of a shortest cycle.)\


In this paper we focus on minimising the function $f$ in
Theorem~\ref{th:Main} and do not calculate $h$ explicitly. In
particular, Theorem~\ref{th:Kt} proves Theorem~\ref{th:Main} with
$f(t)\leq 2^{t-1}+\eps$ for any $\eps>0$ (where the function $h$ also
depends on $\eps$). Note that for Theorem~\ref{th:Main} and all our
results, the proofs can be easily adapted to give polynomial
algorithms that compute the small $K_t$-model.

\begin{figure}[!h]
  \begin{center}
    \includegraphics{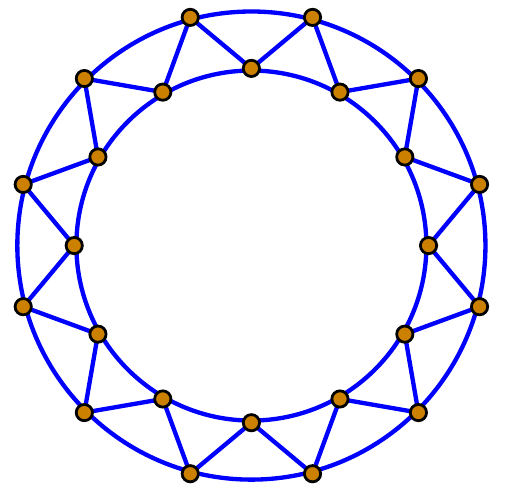}
  \end{center}
  \caption{\label{fig:CycleSquared}$C^2_{24}$}
\end{figure}

For $t\leq 4$, we determine the least possible value of $f(t)$ in
Theorem~\ref{th:Main}. The $t=2$ case is trivial---one edge is a small
$K_2$-minor. To force a small $K_3$-model, average degree $2$ is not
enough, since every $K_3$-model in a large cycle uses every vertex. On
the other hand, we prove that average degree $2+\eps$ forces a cycle
of length $O_{\eps}(\log|G|)$; see Lemma~\ref{lem:cycle}. For $t=4$ we
prove that average degree $4+\eps$ forces a $K_4$-model with
$O_{\eps}(\log|G|)$ vertices; see Theorem~\ref{th:K4}. This result is
also best possible. Consider the square of an even cycle $C_{2n}^2$,
which is a 4-regular graph illustrated in
Figure~\ref{fig:CycleSquared}. If the base cycle is
$(v_1,\dots,v_{2n})$ then $C_{2n}^2-\{v_i,v_{i+1}\}$ is outerplanar
for each $i$. Since outerplanar graphs contain no $K_4$-minor, every
$K_4$-model in $C_{2n}^2$ contains $v_i$ or $v_{i+1}$ for each
$i$, and thus contains at least $n$ vertices.

Motivated by Theorem~\ref{th:Main}, we then consider graphs that
contain $K_3$-models and $K_4$-models of bounded size (not just small
with respect to $|G|$). First, we prove that planar graphs satisfy
this property.  In particular, every planar graph with average degree
at least $2+\eps$ contains a $K_3$-model with $O(\frac{1}{\eps})$
vertices (Theorem~\ref{th:PlanarCycle}). This bound on the average
degree is best possible since a cycle is planar and has average degree
$2$. Similarly, every planar graph with average degree at least
$4+\eps$ contains a $K_4$-model with $O(\frac{1}{\eps})$ vertices
(Theorem~\ref{th:GeneralPlanar}). Again, this bound on the average
degree is best possible since $C^2_{2n}$ is planar and has average
degree $4$. These results generalise for graphs embedded on other
surfaces (Theorems~\ref{th:SurfaceCycle} and \ref{th:SurfaceMinor}).

Finally, we mention three other results in the literature that force a
model of a complete graph of bounded size.

\begin{itemize}

\item \citet{KP88} proved that for every integer $t$ and $\eps>0$,
  every $n$-vertex graph with at least $4^{t^2}n^{1+\eps}$ edges
  contains a subdivision of $K_t$ with at most $\frac{7}{\eps}t^2\log
  t$ vertices; see \citep{Jian-JGT} for recent related results. We
  emphasise that, for fixed $t$, the results in \citep{KP88,Jian-JGT}
  prove that a super-linear lower bound on the number of edges (in
  terms of the number of vertices) forces a $K_t$-minor (in fact, a
  subdivision) of bounded size, whereas Theorem~\ref{th:Main} proves
  that a linear lower bound on the number of edges forces a small
  $K_t$-minor (of size logarithmic in the order of the graph). Also
  note that Theorem~\ref{th:Main} can be proved by adapting the proof
  of \citet{KP88}. As far as we can tell, this method does not give a
  bound better than $f(t)\leq 16^{t}+\eps$ (ignoring lower order
  terms). This bound is inferior to our Theorem~\ref{th:Kt}, which
  proves $f(t)\leq 2^{t-1}+\eps$. Also note that the method of
  \citet{KP88} can be adapted to prove the following result about
  forcing a small subdivision.

\begin{theorem}
  There is a function $h$ such that for every integer $t \geq 2$ and
  real $\eps > 0$, every graph $G$ with average degree at least
  $4^{t^2} + \eps$ contains a subdivision of $K_t$ with at most
  $h(t,\eps)\cdot\log|G|$ division vertices per edge.
\end{theorem}

\item \citet{KO06} proved that every graph with minimum degree at
  least $t$ and girth at least $27$ contains a $K_{t+1}$-subdivision.
  Every graph with average degree at least $2t$ contains a subgraph
  with minimum degree at least $t$. Thus every graph with average
  degree at least $2t$ contains a $K_{t+1}$-subdivision or a
  $K_3$-model with at most $26$ vertices.


\item \citet{KS-GAFA09} proved that for all integers $s'\geq s\geq2$,
  there is a constant $c>0$, such that every $K_{s,s'}$-free graph
  with average degree $r$ contains a minor with average degree at
  least $cr^{1+1/(2s-2)}$.  Applying the result of
  \citet{Kostochka82,Kostochka84} and \citet{Thomason84} mentioned
  above, for every integer $s\geq 2$ there is a constant $c$ such that
  every graph with average degree at least $c(t\sqrt{\log
    t})^{1-1/(2s-1) }$ contains a $K_t$-minor or a $K_{s,s}$-subgraph,
  in which case there is a $K_{s+1}$-model with $2s$ vertices.

\end{itemize}

\section{Definitions and Notations}

See \citep{D05} for undefined graph-theoretic terminology and
notation.  For $S\subseteq V(G)$, let $G[S]$ be the subgraph of $G$
induced by $S$. Let $e(S):=\NumEdges{G[S]}$. For disjoint sets $S,
T\subseteq V(G)$, let $e(S,T)$ be the number of edges between $S$ and
$T$ in $G$.

A \emph{separation} in a graph $G$ is a pair of subgraphs
$\{G_1,G_2\}$, such that $G=G_1\cup G_2$ and $V(G_1)\setminus
V(G_2)\neq\varnothing$ and $V(G_2)\setminus
V(G_1)\neq\varnothing$. The \emph{order} of the separation is
$|V(G_1)\cap V(G_2)|$. A separation of order 1 corresponds to a
cut-vertex $v$, where $V(G_1)\cap V(G_2)=\{v\}$. A separation of order
2 corresponds to a cut-pair $v,w$, where $V(G_1)\cap V(G_2)=\{v,w\}$.

See \citep{MoharThom} for background on graphs embedded in
surfaces. Let $\mathbb{S}_h$ be the orientable surface obtained from
the sphere by adding $h$ handles. The \emph{Euler genus} of
$\mathbb{S}_h$ is $2h$.  Let $\mathbb{N}_c$ be the non-orientable
surface obtained from the sphere by adding $c$ cross-caps. The
\emph{Euler genus} of $\mathbb{N}_c$ is $c$.


An \emph{embedded graph} means a connected graph that is 2-cell
embedded in $\mathbb{S}_h$ or $\mathbb{N}_c$. A \emph{plane graph} is
a planar graph embedded in the plane. Let $F(G)$ denote the set of
faces in an embedded graph $G$.  For a face $f\in F(G)$, let $|f|$ be
the length of the facial walk around $f$.  For a vertex $v$ of $G$,
let $F(G,v)$ be the multiset of faces incident to $v$, where the
multiplicity of a face $f$ in $F(G,v)$ equals the multiplicity of $v$
in the facial walk around $f$. Thus $|F(G,v)|=\deg(v)$.

Euler's formula states that $|G|-\NumEdges{G}+|F(G)|=2-g$ for a
connected graph $G$ embedded in a surface with Euler genus $g$. Note
that $g\leq \NumEdges{G}-|G|+1$ since $|F(G)|\geq1$.  The \emph{Euler
  genus} of a graph $G$ is the minimum Euler genus of a surface in
which $G$ embeds.

We now review some well-known results that will be used implicitly
(see \citep[Section 7.3]{D05}). If a graph $G$ contains no $K_4$-minor
then $\NumEdges{G}\leq 2|G|-3$, and if $|G|\geq2$ then $G$ contains at
least two vertices with degree at most $2$. Hence, if
$\NumEdges{G}>2|G|-3$ then $G$ contains a $K_4$-minor. Similarly, if
$|G|\geq2$ and at most one vertex in $G$ has degree at most $2$, then
$G$ contains a $K_4$-minor.

Throughout this paper, logarithms are binary unless stated otherwise.

\section{Small $K_{3}$-Models and $K_{4}$-Models}

In this section we prove tight bounds on the average degree that
forces a small $K_3$-model or $K_4$-model. The following lemma is at
the heart of many of our results. It is analogous to Lemma~1.1 in
\citep{KP88}

\begin{lemma}
  \label{lem:subgraph}
  There is a function $p$ such that for every two reals $d > d' \geq
  2$, every graph $G$ with average degree at least $d$ contains a
  subgraph with average degree at least $d'$ and diameter at most
  $p(d, d') \cdot \log|G|$.
\end{lemma}

\begin{proof}
  We may assume that every proper subgraph of $G$ has average degree strictly
  less than $d$ (otherwise, simply consider a minimal subgraph with
  that property).  Let
$$\beta := \frac{d}{d'}>1\;\;\;\text{ and }\;\;
p(d, d') := \frac{2}{\log\beta} + 2\enspace.$$ Let $v$ be an arbitrary
vertex of $G$.  Let $B_{k}(v)$ be the subgraph of $G$ induced by the
set of vertices at distance at most $k$ from $v$.  Let $k \geq 1$ be
the minimum integer such that $|B_{k}(v)| < \beta \cdot
|B_{k-1}(v)|$. (There exists such a $k$, since $\beta > 1$ and $G$ is
finite.)\ It follows that $\beta^{k-1} \leq |B_{k-1}(v)| \leq |G|$,
and $B_{k}(v)$ has diameter at most $2k \leq 2(\log_{\beta} |G| + 1)
\leq p(d, d') \cdot \log|G|$.

We now show that $B_{k}(v)$ also has average degree at least $d'$.
Let
\begin{align*}
  A &:=V(B_{k-1}(v)), \\
  B &:=V(B_{k}(v)) \setminus V(B_{k-1}(v)), \\
  C &:= V(G) \setminus (A\cup B)\enspace.
\end{align*}

If $C = \varnothing$, then $B_{k}(v) = G[A\cup B] = G$, and hence
$B_{k}(v)$ has average degree at least $d \geq d'$. Thus, we may
assume that $C\neq\varnothing$.  Let $d''$ be the average degree of
$B_{k}(v)$. Thus,
\begin{equation}
  \label{eq:AB}
  2\big(e(A) + e(B) + e(A,B)\big) = d'' \cdot (|A| + |B|)\enspace.
\end{equation}

Since $C$ is non-empty, $G-A$ is a proper non-empty subgraph of $G$.
By our hypothesis on $G$, this subgraph has average degree strictly
less than $d$; that is,
\begin{equation}
  \label{eq:BC}
  2\big(e(B) + e(C) + e(B,C)\big) < d \cdot (|B| + |C|)\enspace.
\end{equation}

By \eqref{eq:AB} and \eqref{eq:BC} and since $e(A,C)=0$,
\begin{align*}
  2\NumEdges{G} & = 2\big(e(A) + e(B) + e(C) + e(A,B) + e(B,C)\big) \\
  & = d'' (|A|+|B|) + 2e(C) + 2e(B,C) \\
  & < d'' (|A|+|B|) + d (|B|+|C|) - 2e(B) \\
  & \leq d |G| - d |A| + d'' (|A|+|B|)\enspace.
\end{align*}
Thus $d'' (|A|+|B|) > d |A|$ (since $2\NumEdges{G} \geq d\,|G|$).  On
the other hand, by the choice of $k$,
$$\frac{|A|}{|A|+|B|} > \frac{1}{\beta}\enspace.$$
Hence
$$d'' > d \frac{|A|}{|A|+|B|} > \frac{d}{\beta} = d'\enspace,$$
as desired.
\end{proof}




\begin{lemma}
  \label{lem:cycle}
  There is a function $g$ such that for every real $\eps>0$, every
  graph $G$ with average degree at least $2+\eps$ has girth at most
  $g(\eps)\cdot\log|G|$,
\end{lemma}

\begin{proof}
  By Lemma~\ref{lem:subgraph}, $G$ contains a subgraph $G'$ with
  average degree at least $2$ and diameter at most
  $p(2+\eps,2)\cdot\log|G|$.  Let $T$ be a breadth-first search tree
  in $G'$.  Thus $T$ has diameter at most $2p(2+\eps,2)\cdot\log|G|$.
  Since $G'$ has average degree at least 2, $G'$ is not a tree, and
  there is an edge $e\in E(G')\setminus E(T)$. Thus $T$ plus $e$
  contains a cycle of length at most $2 p(2+\eps,2)\cdot\log|G|+1$.
\end{proof}

\begin{theorem}
  \label{th:K4}
  There is a function $h$ such that for every real $\eps > 0$, every
  graph $G$ with average degree at least $4 + \eps$ contains a
  {\md{4}} with at most $h(\eps)\cdot\log|G|$ vertices.
\end{theorem}

\begin{proof}
  By Lemma~\ref{lem:subgraph}, $G$ contains a subgraph $G'$ with
  average degree at least $4 + \frac{\eps}{2}$ and diameter at most
  $p(4 + \eps, 4 + \frac{\eps}{2}) \cdot \log|G|$.  Let $v$ be an
  arbitrary vertex of $G'$.  Let $T$ be a breadth-first search tree
  from $v$ in $G'$.  Let $k$ be the depth of $T$. Thus $k \leq
  p(4+\eps, 4 + \frac{\eps}{2}) \cdot \log|G|$.

  Let $H:= G' - E(T)$. Since $\NumEdges{T}=|G| - 1$, the graph $H$ has
  average degree at least $2 + \frac{\eps}{2}$. By
  Lemma~\ref{lem:cycle}, $H$ contains a cycle $C$ of length at most
  $g(\frac{\eps}{2})\cdot\log|G|$.  We will prove the theorem with
  $h(\eps):=g(\frac{\eps}{2})+3p(4+\eps,4+\frac{\eps}{2})$.

  Observe that $v \notin V(C)$, since $v$ is isolated in $H$.  A
  vertex $w$ of $C$ is said to be \emph{maximal} if, in the tree $T$
  rooted at $v$, no other vertex of $C$ is an ancestor of $w$. Let
  $\dist(x)$ be the distance between $v$ and each vertex $x$ in $T$.

  Consider an edge $xx'$ in $C$ where $x$ is maximal and $x'$ is not.
  Since $T$ is a breadth-first search tree, $\dist(x')\leq\dist(x)+1$.
  Thus, if $x$ is an ancestor of $x'$ then $xx'\in E(T)$, which is a
  contradiction since $xx'\in E(H)$. Hence $x$ is not an ancestor of
  $x'$. Let $y$ be an ancestor of $x'$ in $C$ (which exists since $x'$
  is not maximal). Then $\dist(y)<\dist(x')\leq \dist(x)+1$, implying
  $\dist(y)\leq\dist(x)$. We repeatedly use these facts below.

  First, suppose that there is a unique maximal vertex $x$ in $C$. Let
  $x'$ be a neighbour of $x$ in $C$. Since $x'$ is not maximal, some
  ancestor of $x'$ is in $C$. As proved above, $x$ is not an ancestor
  of $x'$ in $T$, which contradicts the assumption that $x$ is the
  only maximal vertex in $C$.

  Next, suppose there are exactly two maximal vertices $x$ and $y$ in
  $C$.  Let $P$ be an $x$--$y$ path in $C$ that is not the edge $xy$
  (if it exists).  Let $x'$ be the neighbour of $x$ in $P$, and let
  $y'$ be the neighbour of $y$ in $P$. Thus $x'\neq y$ and $y'\neq
  x$. Hence neither $x'$ nor $y'$ are maximal. As proved above, $y$ is
  an ancestor of $x'$ and $\dist(y)\leq\dist(x)$, and $x$ is an
  ancestor of $y'$ and $\dist(x)\leq\dist(y)$. Thus
  $\dist(x)=\dist(y)$. Hence $\dist(x')\leq\dist(y)+1$ and
  $\dist(y')\leq\dist(x)+1$, which implies that $x'y$ and $y'x$ are
  both edges of $T$, and $x' \neq y'$.  Now, the cycle $C$ plus these
  two edges gives a {\md{4}} with $|C|\leq
  g(\frac{\eps}{2})\cdot\log|G|\leq h(\eps)\cdot\log|G|$ vertices.

  Finally, suppose that $C$ contains three maximal vertices $x,y,z$.
  For $w\in \{x,y,z\}$, let $P_{w}$ be the unique $v$--$w$ path in
  $T$.  Then $C \cup P_{x} \cup P_{y} \cup P_{z}$ contains a {\md{4}}
  with at most $|C| + |P_{x} - x| + |P_{y} - y| + |P_{z} - z| \leq |C|
  + 3k \leq h(\eps) \cdot \log|G|$ vertices.
\end{proof}

\section{Small $K_{t}$-Models}

The following theorem establishes our main result
(Theorem~\ref{th:Main}).

\begin{theorem}
  \label{th:WeakKt}
  There is a function $h$ such that for every integer $t \geq 2$ and
  real $\eps > 0$, every graph $G$ with average degree at least $2^t +
  \eps$ contains a {\md{t}} with at most $h(t,\eps) \cdot \log|G|$
  vertices.
\end{theorem}
\begin{proof}
  We prove the following slightly stronger statement: Every graph $G$
  with average degree at least $2^t + \eps$ contains a {\md{t}} with
  at most $h(t, \eps) \cdot \log|G|$ vertices such that each branch
  set of the model contains at least two vertices.

  The proof is by induction on $t$. For $t=2$, let
  $h(t,\eps):=2$. Here we need only assume average degree at least
  $2+\eps$. Some component of $G$ is neither a tree nor a cycle, as
  otherwise $G$ would have average degree at most $2$. It is easily
  seen that this component contains a path on $4$ vertices, yielding a
  {\md{2}} in which each branch set contains two vertices. This model
  has $4 \leq h(t, \eps) \cdot \log|G|$ vertices, as desired. (Observe
  that $|G| \geq 4$, since $G$ contains a vertex with degree at least
  $3$.)\

  Now assume $t \geq 3$ and the claim holds for smaller values of $t$.
  Using Lemma~\ref{lem:subgraph}, let $G'$ be a subgraph of $G$ with
  average degree at least $2^t + \frac{\eps}{2}$ and diameter at most
  $p(2^t+\eps,2^t+\frac{\eps}{2}) \cdot \log|G|$. Let $h(t, \eps) := 2
  + (t-1)\, p(2^t+\eps,2^t+\frac{\eps}{2}) + h(t-1,\frac{\eps}{4})$.

  Choose an arbitrary edge $uv$ of $G'$.  Define the \emph{depth} of a
  vertex $w\in V(G')$ to be the minimum distance in $G'$ between $w$
  and a vertex in $\{u, v\}$.  Note that the depths of the endpoints of
each edge differ by at most $1$. The \emph{depth} of an edge $xy \in
  E(G')$ is the minimum of the depth of $x$ and the depth of $y$.

  Considering edges of $G'$ with even depth on one hand, and with odd
  depth on the other, we obtain two edge-disjoint spanning subgraphs
  of $G'$.  Since $G'$ has average degree at least $2^t +
  \frac{\eps}{2}$, one of these two subgraphs has average degree at
  least $2^{t-1} + \frac{\eps}{4}$.  Let $H$ be a component of this
  subgraph with average degree at least $2^{t-1} + \frac{\eps}{4}$.
  Observe that every edge of $H$ has the same depth $k$ in $G$.

  If $k=0$, then $E(H)$ is precisely the set of edges incident to $u$
  or $v$ (or both). Thus, every vertex in $V(H) \setminus \{u,v\}$ has
  degree at most $2$ in $H$.  Hence $H$ has average degree less than
  $4 < 2^{t-1}+ \frac{\eps}{4}$, a contradiction.  Therefore $k \geq
  1$.

  Now, by induction, $H$ contains a {\md{t-1}} with at most
  $h(t-1,\frac{\eps}{4}) \cdot \log|G'|$ vertices such that each of
  the $t-1$ branch sets $B_{1}, \dots, B_{t-1}$ has at least two
  vertices.  Thus, each $B_i$ contains an edge of $H$. Hence, there is
  a vertex $v_i$ in $B_i$ having depth $k$ in $G'$.  Therefore, there
  is a path $P_{i}$ of length $k$ in $G'$ between $v_{i}$ and some
  vertex in $\{u,v\}$. Let $P_{uv}$ be the trivial path consisting of
  the edge $uv$.  Let
$$
B_{t} := P_{uv} \cup \bigcup_{1 \leq i \leq t-1} (P_{i} -
v_{i})\enspace.$$
The subgraph $B_{t}$ is connected, contains at least two vertices
(namely, $u$ and $v$), and is vertex disjoint from $B_{i}$ for all
$i\in\{1,\dots,t-1\}$.  Moreover, there is an edge between $B_t$ and
each $B_i$, and
\begin{align*}
  \sum_{1 \leq i \leq t} |B_{i}|
  &\leq |B_{t}| + h(t-1,\frac{\eps}{4}) \cdot \log|G'|    \\
  &\leq 2 + \sum_{1 \leq i \leq t-1} |P_{i} - v_{i}| + h(t-1,\frac{\eps}{4}) \cdot \log|G|  \\
  &\leq 2 + (t-1)k + h(t-1,\frac{\eps}{4}) \cdot \log|G|  \\
  &\leq 2 + (t-1)\, p(2^t+\eps,2^t+\frac{\eps}{2}) \cdot \log|G| + h(t-1,\frac{\eps}{4}) \cdot \log|G|  \\
  &\leq h(t,\eps) \cdot \log|G|\enspace.
\end{align*}
Hence, adding $B_{t}$ to our {\md{t-1}} gives the desired {\md{t}} of
$G$.
\end{proof}

Observe that one obstacle to reducing the lower bound on the average
degree in Theorem~\ref{th:WeakKt} is the case $t=3$, which we address
in the following result.

\begin{lemma}
  \label{le:NiceK3} There is a function $h$ such that for every real
  $\eps > 0$, every graph $G$ with average degree at least $4+\eps$
  contains a $K_{3}$-model with at most $h(\eps)\cdot\log|G|$
  vertices, such that each branch set contains at least two vertices.
\end{lemma}

\begin{proof}
  The proof is by induction on $|G| + \NumEdges{G}$.  We may assume
  that no proper subgraph of $G$ has average degree at least $4 +
  \eps$, since otherwise we are done by induction. This implies that
  $G$ is connected.  Note that $|G|\geq6$ since $G$ has average degree
  $> 4$.

  First, suppose that $G$ contains a $K_{4}$ subgraph with vertex set
  $X$.

  {\bf Case 1.} All edges between $X$ and $V(G) \setminus X$ in $G$
  are incident to a common vertex $v\in X$: Let $Y: = X\setminus
  \{v\}$.  Then
$$
2\NumEdges{G-Y} = 2\NumEdges{G} - 12 \geq (4+\eps)|G| - 12 \geq
(4+\eps)|G-Y|\enspace,
$$
implying that $G-Y$ also has average degree at least $4+\eps$, a
contradiction.


{\bf Case 2.} There are two independent edges $uu'$ and $vv'$ between
$X$ and $V(G) \setminus X$ in $G$, where $u, v \in X$: Then $\{u,
u'\}, \{v, v'\}, X \setminus \{u,v\}$ is the desired $K_3$-model.

{\bf Case 3.} Some vertex $w\in V(G) \setminus X$ is adjacent to two
vertices $u, v \in X$: No vertex in $X$ has a neighbour in $V(G)
\setminus (X \cup \{w\})$, as otherwise Case 2 would apply.  Since $G$
is connected and $|G| \geq 6$, it follows that $w$ has a neighbour
$w'$ outside $X$. Let $x, y$ be the two vertices in $X \setminus \{u,
v\}$.  Then $\{w, w'\}, \{u, x\}, \{v, y\}$ is the desired
$K_3$-model.

This concludes the case in which $G$ contains a $K_{4}$ subgraph.
Now, assume that $G$ is $K_{4}$-free. By Theorem~\ref{th:K4}, $G$
contains a $K_4$-model $B_1,\dots,B_4$ with at most
$h(\eps)\cdot\log|G|$ vertices. Without loss of generality,
$|B_1|\geq|B_2|\geq|B_3|\geq|B_4|$ and $|B_1|\geq2$.

{\bf Case 1.} $|B_2|\geq2$: Then $B_1,B_2,B_3\cup B_4$ is the desired
$K_3$-model. Now assume that $B_i=\{x_i\}$ for all $i\in\{2,3,4\}$.

{\bf Case 2.} Some $x_i$ is adjacent to some vertex $w$ not in
$B_1\cup B_2\cup B_3\cup B_4$: If $i=2$ then $\{x_2,w\},B_1,B_3\cup
B_4$ is the desired $K_3$-model. Similarly for $i\in\{3,4\}$.


{\bf Case 3.} $|B_1|\geq 3$.  Then there are two independent edges in
$G$ between $B_{1}$ and $\{x_{2}, x_{3}, x_{4}\}$, say $ux_{2}$ and
$vx_{3}$ with $u,v \in B_{1}$ (otherwise, there would be a $K_{4}$
subgraph).  There is a vertex $w\in B_{1} \setminus \{u,v\}$ adjacent
to at least one of $u,v$, say $u$. Let $C$ be the vertex set of the
component of $G[B_{1}] - \{u,w\}$ containing $v$.  Then $\{u, w\}, C
\cup \{x_{3}\}, \{x_{2}, x_{4}\}$ is the desired $K_3$-model.

{\bf Case 4.} $B_1 = \{u, v\}$.  As in the previous cases, there are
two independent edges in $G$ between $\{u,v\}$ and $\{x_{2}, x_{3},
x_{4}\}$, say $ux_{2}$ and $vx_{3}$.  At least one of $u,v$, say $u$,
is adjacent to some vertex $w$ outside $\{u, v, x_{2}, x_{3},
x_{4}\}$, because $G$ is connected with at least $6$ vertices, and
none of $x_{2}$, $x_{3}$, $x_{4}$ has a neighbour outside $\{u, v,
x_{2}, x_{3}, x_{4}\}$.  Then $\{u, w\}, \{v, x_{3}\}, \{x_{2},
x_{4}\}$ is the desired $K_3$-model.
\end{proof}

Note that average degree greater than $4$ is required in
Lemma~\ref{le:NiceK3} because of the disjoint union of $K_5$'s.
Lemma~\ref{le:NiceK3} enables the following improvement to
Theorem~\ref{th:WeakKt}.

\begin{theorem}
  \label{th:Kt}
  There is a function $h$ such that for every integer $t \geq 2$ and
  real $\eps > 0$, every graph $G$ with average degree at least
  $2^{t-1} + \eps$ contains a {\md{t}} with at most $h(t,
  \eps)\cdot\log|G|$ vertices.
\end{theorem}
\begin{proof}
  As before, we prove the following stronger statement: Every graph
  $G$ with average degree at least $2^{t-1} + \eps$ contains a
  {\md{t}} with at most $h(t, \eps) \cdot \log|G|$ vertices such that
  each branch set of the model contains at least two vertices.

  The proof is by induction on $t$. The $t=2$ case is handled in the
  proof of Theorem~\ref{th:WeakKt}. Lemma~\ref{le:NiceK3} implies the
  $t=3$ case.  Now assume $t \geq 4$ and the claim holds for smaller
  values of $t$.  The proof proceeds as in the proof of
  Theorem~\ref{th:WeakKt}. We obtain a subgraph $G'$ of $G$ with
  average degree at least $2^{t-1} + \frac{\eps}{2}$ and diameter at
  most $p(2^{t-1}+\eps,2^{t-1}+\frac{\eps}{2}) \cdot \log|G|$. Choose
  an edge $uv$ of $G'$ and define the depth of edges with respect to
  $uv$. We obtain a connected subgraph $H$ with average degree at
  least $2^{t-2} + \frac{\eps}{4}$, such that every edge of $H$ has
  the same depth $k$. If $k=0$, then $E(H)$ is precisely the set of
  edges incident to $u$ or $v$ (or both), implying $H$ has average
  degree less than $4 < 2^{t-2}+ \frac{\eps}{4}$. Now assume
  $k\geq1$. The remainder of the proof is the same as that of
  Theorem~\ref{th:WeakKt}.
\end{proof}


\citet{Thomassen-JCTB83} first observed that high girth (and minimum
degree 3) forces a large complete graph as a minor; see \citep{KO03}
for the best known bounds. We now show that high girth (and minimum
degree 3) forces a \emph{small} model of a large complete graph.

\begin{theorem} 
  Let $k$ be a positive integer. Let $G$ be a graph with girth at
  least $8k+3$ and minimum degree $r\geq3$. Let $t$ be an integer such
  that $r(r-1)^k\geq 2^{t-1}+1$. Then $G$ contains a $K_t$-model with
  at most $h'(k,r)\cdot\log|G|$ vertices, for some function $h'$.
\end{theorem}

\begin{proof}
  \citet{Mader-Comb98} proved that $G$ contains a minor $H$ of minimum
  degree at least $r(r-1)^k$, such that each branch set has radius at
  most $2k$; see \citep[Lemma~7.2.3]{D05}. Let
  $V(H)=\{b_1,\dots,b_{|H|}\}$, and let $B_1,\dots,B_{|H|}$ be the
  corresponding branch sets in $G$. Let $r_i$ be a centre of
  $B_i$. For each vertex $v$ in $B_i$, let $P_{i,v}$ be a path between
  $r_i$ and $v$ in $B_i$ of length at most $2k$.

  By Theorem~\ref{th:Kt}, $H$ contains a $K_t$-model with at most
  $h(t)\cdot\log|H|$ vertices. Let $C_1,\dots,C_t$ be the
  corresponding branch sets. Say $C_i$ has $n_i$ vertices. Thus
  $\sum_{i=1}^tn_i\leq h(t)\cdot\log|H|$. We now construct a
  $K_t$-model $X_1,\dots,X_t$ in $G$.

  For $i\in\{1,\dots,t\}$, let $T_i$ be a spanning tree of $C_i$.
  Each edge $b_jb_\ell$ of $T_i$ corresponds to an edge $vw$ of $G$,
  for some $v$ in $B_j$ and $w$ in $B_\ell$.  Add to $X_i$ the
  $r_ir_j$-path $P_{j,v}\cup\{vw\}\cup P_{\ell,w}$.  This path has at
  most $4k+2$ vertices. Thus $X_i$ is a connected subgraph of $G$ with
  at most $(4k+2)(n_i-1)$ vertices (since $T_i$ has $n_i-1$ edges).

  For distinct $i,i'\in\{1,\dots,t\}$ there is an edge between $C_i$
  and $C_{i'}$ in $H$. This edge corresponds to an edge $vw$ of $G$,
  where $v$ is in some branch set $B_j$ in $C_i$, and $w$ is in some
  branch set $B_{j'}$ in $C_{i'}$.  Add the path $P_{j,v}$ to $X_i$,
  and add the path $P_{j',w}$ to $X_{i'}$. Thus $v$ in $X_i$ is
  adjacent to $w$ in $X_j$.

  Hence $X_1,\dots,X_t$ is a $K_t$-model in $G$ with at most
  $\sum_{i=1}^t(4k+2)(n_i-1)\leq (4k+2)\cdot h(t)\cdot\log|H|$
  vertices from the first step of the construction, and at most
  $\binom{t}{2}(4k+2)$ vertices from the second step.  Since $t$ is
  bounded by a function of $r$ and $k$, there are at most
  $h'(k,r)\cdot\log|G|$ vertices in total, for some function $h'$.
\end{proof}

\begin{corollary}
  Let $k$ be a positive integer.  Let $G$ be a graph with girth at
  least $8k+3$ and minimum degree at least $3$. Then $G$ contains a
  $K_k$-model with at most $h(k)\cdot\log|G|$ vertices, for some
  function $h$.\qed
\end{corollary}

\section{Planar Graphs}
\label{sec:Planar}

In this section we prove that sufficiently dense planar graphs have
$K_3$-models and $K_4$-models of bounded size. We start with the $K_3$
case.

\begin{theorem}
  \label{th:PlanarCycle}
  Let $\eps\in(0,4)$. Every planar graph $G$ with average degree at
  least $2+\eps$ has girth at most $1+\CeilFrac{4}{\eps}$.
\end{theorem}

\begin{proof}
  Let $H$ be a connected component of $G$ with average degree at least
  $2+\eps$. Thus $H$ is not a tree. Say $H$ has $n$ vertices and $m$
  edges. Fix an embedding of $H$ in the plane with $r$ faces. Let
  $\ell$ be the minimum length of a facial walk. Thus $\ell\geq3$ and
  $2m\geq r\ell=(2+m-n)\ell$, implying
$$n-2\geq m(1-\tfrac2\ell)\geq \tfrac12(2+\eps)n(1-\tfrac2\ell)> \tfrac12(2+\eps)(n-2)(1-\tfrac2\ell)\enspace.$$
It follows that $\ell<2+\frac{4}{\eps}$.  Since $\ell$ is an integer,
$\ell\leq 1+\CeilFrac{4}{\eps}$.  Since $H$ is not a tree, every
facial walk contains a cycle.  Thus $H$ and $G$ have girth at most
$1+\CeilFrac{4}{\eps}$.
\end{proof}


To prove our results for $K_4$-models in embedded graphs, the notion
of visibility will be useful (and of independent interest). Distinct
vertices $v$ and $w$ in an embedded graph are \emph{visible} if $v$
and $w$ appear on a common face; we say $v$ \emph{sees} $w$.

\begin{lemma}
  \label{lem:FindMinor}
  Let $v$ be a vertex of a plane graph $G$, such that $\deg(v)\geq3$,
  $v$ is not a cut-vertex, and $v$ is in no cut-pair. Then $v$ and the
  vertices seen by $v$ induce a subgraph containing a $K_4$-minor.
\end{lemma}

\begin{proof}
  We may assume that $G$ is connected.  Since $v$ is not a cut-vertex,
  $G-v$ is connected.  Let $f$ be the face of $G-v$ that contains $v$
  in its interior.  Let $F$ be the facial walk around $f$.  Suppose
  that $F$ is not a simple cycle.  Then $F$ has a repeated vertex $w$.
  Say $(a,w,b,\dots,c,w,d)$ is a subwalk of $F$.  Then there is a
  Jordan curve $C$ from $v$ to $w$, arriving at $w$ between the edges
  $wa$ and $wb$, then leaving $w$ from between the edges $wc$ and
  $wd$, and back to $v$.  Thus $C$ contains $b$ in its interior and
  $a$ in its exterior.  Hence $v,w$ is a cut-pair.  This contradiction
  proves that $F$ is a simple cycle.  Hence $v$ and the vertices seen
  by $v$ induce a subdivided wheel with $\deg(v)$ spokes.  Since
  $\deg(v)\geq 3$ this subgraph contains a subdivision of $K_4$.
\end{proof}

Recall that $F(G,v)$ is the multiset of faces incident to a vertex $v$
in an embedded graph $G$, where the multiplicity of a face $f$ in
$F(G,v)$ equals the multiplicity of $v$ in the facial walk around $f$.

\begin{lemma}
  \label{lem:Sees}
  Each vertex $v$ in an embedded graph $G$ sees at most
  $$\sum_{f\in F(G,v)}(|f|-2)$$ other vertices.
\end{lemma}

\begin{proof}
  The vertex $v$ only sees the vertices in the faces in $F(G,v)$. Each
  $f\in F(G,v)$ contributes at most $|f|-1$ vertices distinct from
  $v$. Moreover, each neighbour of $v$ is counted at least twice. Thus
  $v$ sees at most $\sum_{f\in F(G,v)}(|f|-1)-\deg(v)$ other vertices,
  which equals $\sum_{f\in F(G,v)}(|f|-2)$.
\end{proof}

The 4-regular planar graph $C^2_{2n}$ has an embedding in the plane,
in which each vertex sees $n+1$ other vertices; see
Figure~\ref{fig:CycleSquared}. On the other hand, we now show that
every plane graph with minimum degree 5 has a vertex that sees a
bounded number of vertices.

\begin{lemma}
  \label{MinimumDegreeFiveSees}
  Every plane graph $G$ with minimum degree $5$ has a vertex that sees
  at most $7$ other vertices.
\end{lemma}

\begin{proof}
  For each vertex $v$ of $G$, associate a charge of
$$2-\deg(v)+\sum_{f\in F(G,v)}\frac{2}{|f|}\enspace.$$
By Euler's formula, the total charge is
$2|G|-2\NumEdges{G}+2|F(G)|=4$.  Thus some vertex $v$ has positive
charge. That is,
\begin{align*}
  2\sum_{f\in F(G,v)}\frac{1}{|f|}>\deg(v)-2 \enspace.
\end{align*}
Now $\frac{1}{|f|}\leq\frac{1}{3}$.  Thus
$\frac{2}{3}\deg(v)>\deg(v)-2$, implying $\deg(v)<6$ and $\deg(v)=5$.
If some facial walk containing $v$ has length at least $6$, then
$$3=2\left(\frac{4}{3}+\frac{1}{6} \right) \geq 2 \sum_{f\in F(G,v)} \frac{1 }{|f|} > 3\enspace,$$
which is a contradiction. Hence each facial walk containing $v$ has
length at most $5$.  If two facial walks containing $v$ have length at
least $4$, then
$$3=2\left(\frac{3}{3}+\frac{2}{4} \right) \geq 2 \sum_{f\in F(G,v)} \frac{1 }{|f|} > 3\enspace,$$
which is a contradiction. Thus no two facial walks containing $v$ each
have length at least $4$. Hence all the facial walks containing $v$
are triangles, except for one, which has length at most $5$.  Thus $v$
sees at most $7$ vertices.
\end{proof}

\begin{figure}[!ht]
  \begin{center}
    \includegraphics{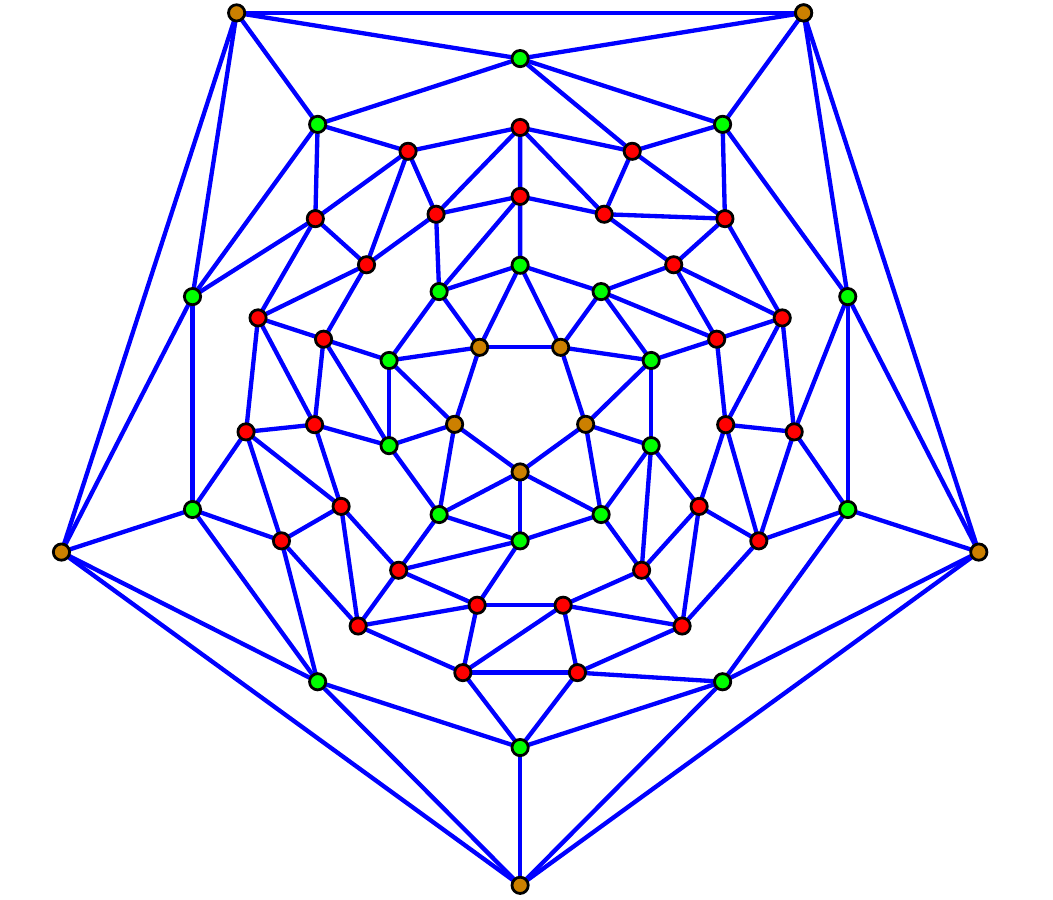}
  \end{center}
  \caption{\label{fig:60}The snub dodecahedron.}
\end{figure}

The bound in Lemma~\ref{MinimumDegreeFiveSees} is tight since there is
a 5-regular planar graph with triangular and pentagonal faces, where
each vertex is incident to exactly one pentagonal face (implying that
each vertex sees exactly 7 vertices).  The corresponding polyhedron is
called the \emph{snub dodecahedron}; see Figure~\ref{fig:60}.

Lemmas~\ref{lem:FindMinor} and \ref{MinimumDegreeFiveSees} imply:

\begin{theorem}
  \label{th:3ConnMinDeg5}
  Every 3-connected planar graph with minimum degree $5$ contains a
  $K_4$-model with at most $8$ vertices.
\end{theorem}

Theorem~\ref{th:3ConnMinDeg5} is best possible since it is easily seen
that every $K_4$-model in the snub dodecahedron contains at least 8
vertices. Also note that no result like Theorem~\ref{th:3ConnMinDeg5}
holds for planar graphs with minimum degree 4 since every $K_4$-model
in the 4-regular planar graph $C_{2n}^2$ has at least $n$ vertices.



We now generalise Lemma~\ref{MinimumDegreeFiveSees} for graphs with
average degree greater than $4$.

\begin{lemma}
  \label{AverageDegreeFourPlusSees}
  Let $\eps\in(0,2)$.  Every plane graph $G$ with minimum degree at
  least $3$ and average degree at least $4+\eps$ has a vertex $v$ that
  sees at most $1+\ceil{\frac{8}{\eps}}$ other vertices.
\end{lemma}

\begin{proof}
  For each vertex $v$ of $G$, associate a charge of
$$(8+2\eps)-\Brac{8+3\eps}\deg(v)+ \Brac{24+6\eps}\sum_{f\in F(G,v)}\frac{ 1}{|f|}\enspace. $$
By Euler's formula, the total charge is
\begin{align*}
  \;&(8+2\eps)|G|-\Brac{16+6\eps}\NumEdges{G}+ \Brac{24+6\eps}|F(G)|\\
  =\;&(8+2\eps)|G|-\Brac{16+6\eps}\NumEdges{G}+ \Brac{24+6\eps}(\NumEdges{G}-|G|+2)\\
  =\;&4(2\NumEdges{G}-(4+\eps)|G|)+ 2\Brac{24+6\eps}\\
  \geq\;& 2\Brac{24+6\eps} \enspace.
\end{align*}
Thus some vertex $v$ has positive charge.  That is,
\begin{align*}
  (24+6\eps) \sum_{f\in F(G,v)}\frac{ 1}{|f|} > ( 8+3\eps) \deg(v)
  -(8+2\eps) \enspace.
\end{align*}
That is,
\begin{align*}
  \sum_{f\in F(G,v)}\frac{ 1}{|f|} >
  \Brac{\frac{1}{3}+\frac{1}{\alpha}} \deg(v) -\frac{1}{3} \enspace,
\end{align*}
where $\alpha:=6+\frac{24}{\eps}$.  We have proved that $\deg(v)$ and
the lengths of the facial walks incident to $v$ satisfy
Lemma~\ref{lem:technical} in Appendix~\ref{app}.  Thus $$\sum_{f\in
  F(G,v)}(|f|-2) \leq \left\lceil \frac{\alpha}{3} \right\rceil -1 =
1+\Ceil{\frac{8}{\eps}} \enspace.$$ The result follows from
Lemma~\ref{lem:Sees}.
\end{proof}


Lemmas~\ref{AverageDegreeFourPlusSees} and \ref{lem:FindMinor} imply:

\begin{theorem}
  \label{th:3ConnAverageDegree4+}
  Let $\eps\in(0,2)$.  Every 3-connected planar graph $G$ with average
  degree at least $4+\eps$ contains a $K_4$-model with at most
  $2+\Ceil{\frac{8}{\eps}}$ vertices.
\end{theorem}

We now prove that the 3-connectivity assumption in
Theorem~\ref{th:3ConnAverageDegree4+} can be dropped, at the expense
of a slightly weaker bound on the size of the $K_4$-model.

\begin{theorem}
  \label{th:GeneralPlanar}
  Let $\eps\in(0,2)$.  Every planar graph $G$ with average degree at
  least $4+\eps$ contains a $K_4$-model with at most
  $\ceil{\frac{8}{\eps}}+\ceil{\frac{2}{\eps}}$ vertices. Moreover,
  this bound is within a constant factor of being optimal.
\end{theorem}

\begin{figure}[!ht]
  \begin{center}
    \includegraphics{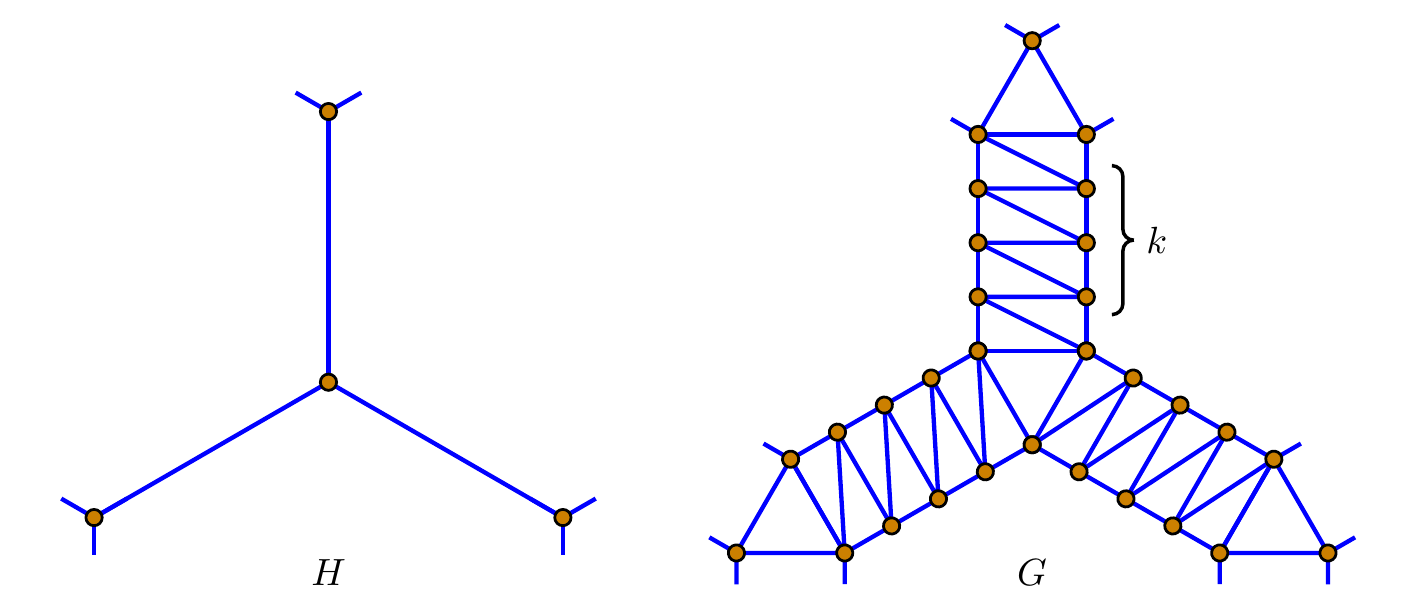}
  \end{center}
  \caption{\label{fig:Construction}Construction of $G$.}
\end{figure}

\begin{proof}
  If $G$ has at most $2+\ceil{\frac{2}{\eps}}$ vertices, then we are
  done since $m>2n$ implies $G$ contains a $K_4$-model, which
  necessarily has at most
  $2+\ceil{\frac{2}{\eps}}<\ceil{\frac{8}{\eps}}+\ceil{\frac{2}{\eps}}$
  vertices.

  We now proceed by induction on $n$ with the following hypothesis:
  Let $G$ be a planar graph with $n\geq 2+\ceil{\frac{2}{\eps}}$
  vertices and $m$ edges, such that
  \begin{align}
    \label{eqn:GeneralPlanar}
    2m>(4+\eps)(n-2)\enspace.
  \end{align}
  Then $G$ contains a $K_4$-model with at most
  $\ceil{\frac{8}{\eps}}+\ceil{\frac{2}{\eps}}$ vertices.

  This will imply the theorem since $2m\geq (4+\eps)n>(4+\eps)(n-2)$.

  Suppose that $n\leq\ceil{\frac{8}{\eps}}+\ceil{\frac{2}{\eps}}$.
  Since $n\geq2+\frac{2}{\eps}$,
$$2m>(4+\eps)(n-2)=4n-8+\eps(n-2) \geq 4n-6\enspace. $$
Thus $m>2n-3$, implying $G$ contains a $K_4$-model, which has at most
$n\leq\ceil{\frac{8}{\eps}}+\ceil{\frac{2}{\eps}}$ vertices. Now
assume that $n\geq\ceil{\frac{8}{\eps}}+\ceil{\frac{2}{\eps}}+1$.

Suppose that $\deg(v)\leq2$ for some vertex $v$. Thus $G-v$ satisfies
\eqref{eqn:GeneralPlanar} since
$$2\NumEdges{G-v}=2(m-\deg(v))>(4+\eps)(n-2)-4>(4+\eps)(n-3)\enspace.$$
Now $n-1\geq
\ceil{\frac{8}{\eps}}+\ceil{\frac{2}{\eps}}>2+\ceil{\frac{2}{\eps}}$. Thus,
by induction, $G-v$ and hence $G$ contains the desired
$K_4$-minor. Now assume that $\deg(v)\geq3$ for every vertex $v$.

Suppose that $G$ contains a separation $\{G_1,G_2\}$ of order at most
$2$. Let $S:=V(G_1\cap G_2)$. Say each $G_i$ has $n_i$ vertices and
$m_i$ edges. Thus $n_1+n_2\leq n+2$ and $m_1+m_2\geq m$.
Equation~\eqref{eqn:GeneralPlanar} is satisfied for $G_1$ or $G_2$, as
otherwise
$$(4+\eps)(n-2) < 2m \leq  2m_1+2m_2\leq
(4+\eps)(n_1+n_2-4)\leq (4+\eps)(n-2)\enspace.$$ Without loss of
generality, $G_1$ satisfies \eqref{eqn:GeneralPlanar}.  Thus we are
done by induction if $n_1\geq2+\ceil{\frac{2}{\eps}}$.  Now assume
that $n_1\leq 1+\ceil{\frac{2}{\eps}}$.  Also assume that $m_1\leq
2n_1-3$, as otherwise $G_1$ contains a $K_4$-model, which has at most
$n_1\leq 1+\ceil{\frac{2}{\eps}}$ vertices.

Suppose that $S=\{v\}$ for some cut-vertex $v$.  Since every vertex in
$G$ has degree at least 3, every vertex in $G_1$, except $v$, has
degree at least $3$ in $G_1$.  Since $n_1\geq2$, $G_1$ contains a
$K_4$-model, which has at most $n_1\leq1+\ceil{\frac{2}{\eps}}$
vertices. Now assume that $G$ is 2-connected.

Suppose that $S=\{v,w\}$ for some adjacent cut-pair $v,w$.  Thus
$n_1+n_2=n+2$ and $m=m_1+m_2-1$ and
\begin{align*}
  2m_2 =2m+2-2m_1 >(4+\eps)(n-2)+2-2(2n_1-3)
  =\;&(4+\eps)(n_1+n_2-4)-4n_1+8\\
  =\;&(4+\eps)(n_2-4)+\eps n_1+8\\
  \geq\;& (4+\eps)(n_2-4)+2(4+\eps)\\
  =\;& (4+\eps)(n_2-2)\enspace.
\end{align*}
That is, $G_2$ satisfies \eqref{eqn:GeneralPlanar}.
Also, $$n_2=n-n_1+2
\geq\Brac{\Ceil{\frac{8}{\eps}}+\Ceil{\frac{2}{\eps}}} +1 -
\Brac{1+\Ceil{\frac{2}{\eps}}} +2 =2+\Ceil{\frac{8}{\eps}}
>2+\Ceil{\frac{2}{\eps}}\enspace.$$ Hence, by induction $G_2$ and thus
$G$ contains the desired $K_4$-model.  Now assume that every cut-pair
of vertices are not adjacent.

Suppose that $S=\{v,w\}$ for some non-adjacent cut-pair $v,w$ and
$m_1\leq 2n_1-4$: Thus $n_1+n_2=n+2$ and $m_1+m_2=m$ and
\begin{align*}
  2m_2 =2m-2m_1 >(4+\eps)(n-2)-2(2n_1-4)
  =\;&(4+\eps)(n_1+n_2-4)-4n_1+8\\
  =\;&(4+\eps)(n_2-4)+\eps n_1+8\\
  \geq\;& (4+\eps)(n_2-4)+2\eps+8\\
  =\;& (4+\eps)(n_2-2)\enspace.
\end{align*}
That is, $G_2$ satisfies \eqref{eqn:GeneralPlanar}.  As proved above,
$n_2 >2+\ceil{\frac{2}{\eps}}$.  Hence, by induction $G_2$ and thus
$G$ contains the desired $K_4$-model.  Now assume that for every
cut-pair $v,w$ we have $vw\not\in E(G)$, and if $\{G_1,G_2\}$ is the
corresponding separation with $G_1$ satisfying
\eqref{eqn:GeneralPlanar}, then $m_1=2n_1-3$ and $n_1\leq
1+\ceil{\frac{2}{\eps}}$.

Fix an embedding of $G$. By Lemma~\ref{AverageDegreeFourPlusSees},
there is a vertex $v$ in $G$ that sees at most
$1+\Ceil{\frac{8}{\eps}}$ other vertices.  If $v$ is in no cut-pair
then by Lemma~\ref{lem:FindMinor} and since $G$ is 2-connected, $v$
plus the vertices seen by $v$ induce a subgraph that contains a
$K_4$-model, which has at most
$2+\Ceil{\frac{8}{\eps}}\leq\Ceil{\frac{8}{\eps}}+\Ceil{\frac{2}{\eps}}$
vertices.  Now assume that $v,w$ is a cut-pair. Thus $vw\not\in E(G)$,
and if $\{G_1,G_2\}$ is the corresponding separation, then
$m_1=2n_1-3$ and $n_1\leq 1+\ceil{\frac{2}{\eps}}$.  Since $v,w$ is a
cut-pair, there is a $vw$-path $P$ contained in $G_2$, such that $P$
is contained in a single face of $G$.  Every vertex in $P$ is seen by
$v$, and $v$ sees at least 2 vertices in $G_1-w$. Thus $P$ has at most
$\Ceil{\frac{8}{\eps}}-2$ internal vertices.  Let $H$ be the minor of
$G$ obtained by contracting $P$ into the edge $vw$, and deleting all
the other vertices in $G_2$. Thus $H$ has $n_1$ vertices and $2n_1-2$
edges. Hence $H$ contains a $K_4$-minor.  The corresponding
$K_4$-model in $G$ is contained in $G_1\cup P$, and thus has at most
$(1+\ceil{\frac{2}{\eps}})+(\ceil{\frac{8}{\eps}}-2)
<\ceil{\frac{2}{\eps}}+\ceil{\frac{8}{\eps}}$ vertices.

We now prove the lower bound. Assume that $\eps\in(0,1]$ and
$k:=\frac{1}{\eps}-1$ is a non-negative integer. Let $H$ be a cubic
plane graph in which the length of every facial walk is at least $5$
(for example, the dual of a minimum degree $5$ plane
triangulation). Say $H$ has $p$ vertices. Let $G$ be the plane graph
obtained by replacing each vertex of $H$ by a triangle, and replacing
each edge of $H$ by $2k$ vertices, as shown in
Figure~\ref{fig:Construction}. Thus $G$ has $3p$ vertices with degree
$5$ and $3kp$ vertices with degree $4$.  Thus
$|G|=3p+3pk=\frac{3p}{\eps}$ and $2\NumEdges{G}=3p\cdot 5+ 3pk \cdot
4=4|G|+3p=(4+\eps)|G|$. Thus $G$ has average degree $4+\eps$.  Every
$K_4$-model in $G$ includes a cycle that surrounds a `big' face with
more than $5k$ vertices. Thus every $K_4$-model has more than
$5k=\frac{5}{\eps}-5$ vertices. Similar constructions are possible for
$\eps>1$ starting with a 4- or 5-regular planar graph.
\end{proof}

\section{Higher Genus Surfaces}

We now extend our results from Section~\ref{sec:Planar} for graphs
embedded on other surfaces.

\begin{lemma}
  \label{le:SurfaceFace}
  Let $\eps>0$. Let $G$ be a graph with average degree at least
  $2+\eps$.  Suppose that $G$ is embedded in a surface with Euler
  genus at most $g$. Then some facial walk has length at most
  $(\frac{4}{\eps}+2)(g+1)$. Moreover, this bound is tight up to lower
  order terms.
\end{lemma}

\begin{figure}[!ht]
  \begin{center}
    \includegraphics{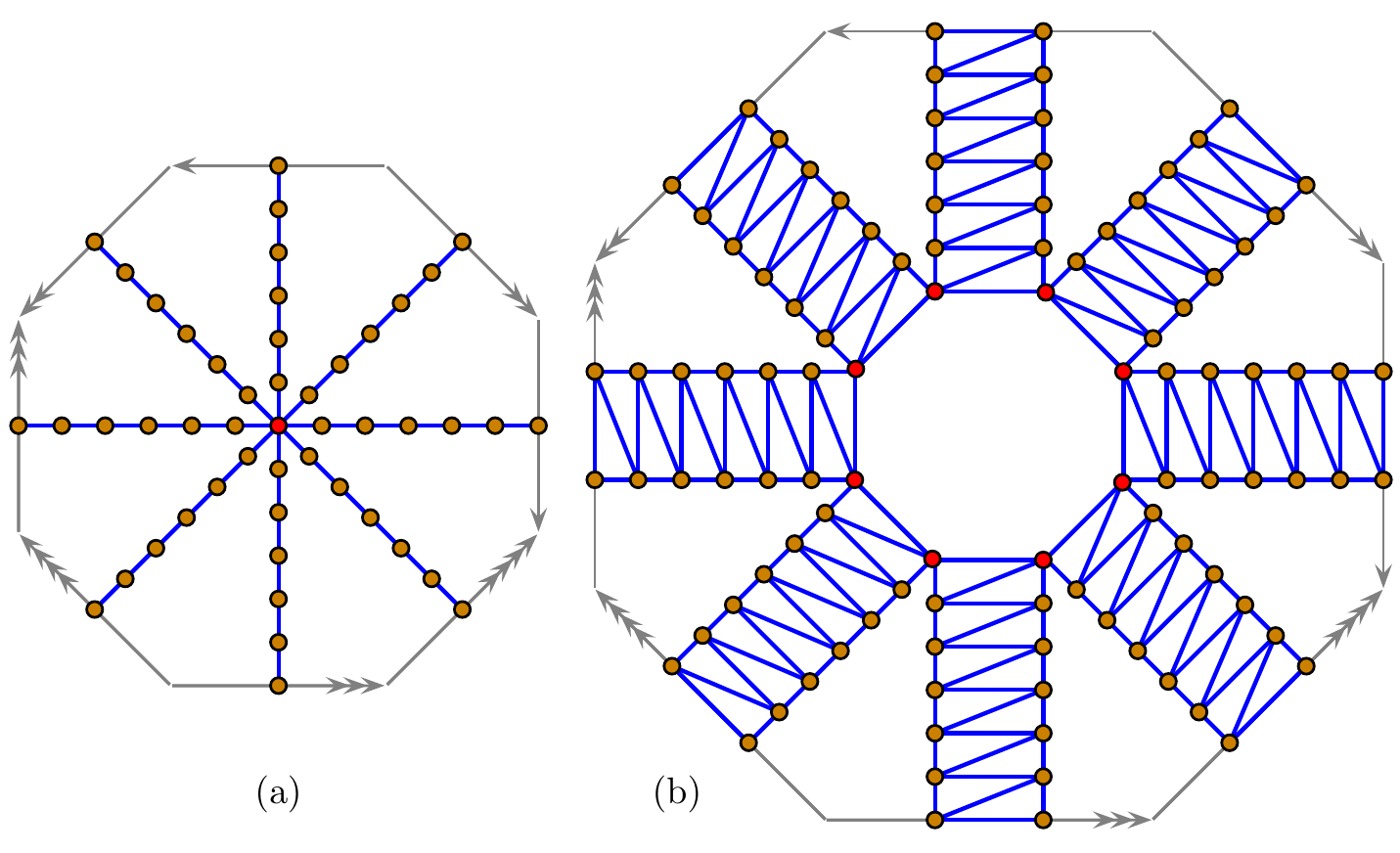}
  \end{center}
  \caption{\label{fig:Surface} Graphs embedded in $\mathbb{S}_2$: (a)
    average degree $2+\eps$ and one face, and (b) average degree
    $4+\eps$ and every vertex on one face.}
\end{figure}

\begin{proof}
  Say $G$ has $n$ vertices, $m$ edges, and $r$ faces.  Let $\ell$ be
  the minimum length of a facial walk. Thus $2m\geq r\ell$. By Euler's
  formula, $n-m+r=2-g$.  Hence
  \begin{align*}
    (2+\eps)n & \leq 2m\\
    (2+\eps)(2-g) & = (2+\eps) (n-m+r) \\
    \frac{\eps}{2}( r\ell) & \leq \frac{\eps}{2} (2m)\enspace.
  \end{align*}
  Summing gives $\frac{\eps}{2}(r\ell) \leq (2+\eps)(g+r-2) $.  Since
  $r\geq1$,
$$\ell \leq \frac{2}{\eps r}\Brac{2+\eps}\Brac{g+r-2} = \Brac{\frac{4}{\eps}+2}\Brac{\frac{g}{r}
  + \frac{r-2}{r}} < \Brac{\frac{4}{\eps} +2} \Brac{g+1} \enspace.$$
Hence some facial walk has length at most $(\frac{4}{\eps}+2)(g+1)$.

Now we prove the lower bound. Assume that $g=2h\geq2$ is a positive
even integer, and that $0<\eps\leq 1 - \frac{3}{2g+1}$. Let
$k:=\floor{\frac{2}{\eps}-\frac{2}{\eps g}-\frac{1}{g}}$. Thus $k\geq
2$. Let $G$ be the graph consisting of $g$ cycles of length $k+1$ with
exactly one vertex in common. Thus
\begin{align*}
  2\NumEdges{G}=2g(k+1)= 2gk + 2 + \eps+ \eps
  g\Brac{\frac{2}{\eps}-\frac{2}{\eps g}-\frac{1}{g}}
  &\geq 2gk + 2 + \eps+  \eps gk \\
  &=  (2+\eps)(gk+1)\\
  &= (2+\eps)|G| \enspace.
\end{align*}
Hence $G$ has average degree at least $2+\eps$.  As illustrated in
Figure~\ref{fig:Surface}(a), $G$ has an embedding in $\mathbb{S}_h$
(which has Euler genus $2h=g$) with exactly one face. Thus every
facial walk in $G$ has length
$2\NumEdges{G}=2g(k+1)>2g(\frac{2}{\eps}-\frac{2}{\eps
  g}-\frac{1}{g})\geq\frac{4(g-1)}{\eps}-2$.
\end{proof}

\begin{theorem}
  \label{th:SurfaceCycle}
  There is a function $h$, such that for every real $\eps>0$, every
  graph $G$ with average degree at least $2+\eps$ and Euler genus $g$
  has girth at most $h(\eps)\cdot \log(g+2)$. Moreover, for fixed
  $\eps$, this bound is within a constant factor of being optimal.
\end{theorem}

\begin{proof}
  Say $G$ has $n$ vertices and $m$ edges. We may assume that every
  proper subgraph of $G$ has average degree strictly less than
  $2+\eps$.  This implies that $G$ has minimum degree at least 2. Fix
  an embedding of $G$ with Euler genus $g$. Let $\ell$ be the minimum
  length of a facial walk. By Euler's formula, there are $m-n+2-g$
  faces. Thus $2m\geq (m-n+2-g)\ell$, implying $\ell(n+g-2)\geq
  m(\ell-2)\geq\half(2+\eps)(\ell-2)n$. Thus $\ell(n+g-2)\geq
  \half(2+\eps)(\ell-2)n$, implying $\ell(g-2)\geq(\frac{\eps}{2}
  (\ell-2) -2)n$. First suppose that $\ell<6+\tfrac{12}{\eps}$. Since
  $G$ has no degree-1 vertices, every facial walk contains a
  cycle. Thus $G$ has girth at most $6+\tfrac{12}{\eps}$, which is at
  most $h(\eps)\cdot \log(g+2)$ for some function $h$. Now assume that
  $\ell\geq 6+\tfrac{12}{\eps}$, which implies that
  $\ell(g-2)\geq(\frac{\eps}{2} (\ell-2) -2)n\geq\frac{\eps}{3} \ell
  n$. Thus $n\leq \frac{3}{\eps}(g-2)$.  By Lemma~\ref{lem:cycle}, the
  girth of $G$ is at most $g(\eps)\cdot\log n\leq
  g(\eps)\cdot\log(\frac{3}{\eps}(g-2))$, which is at most
  $h(\eps)\cdot \log(g+2)$ for some function $h$.

  Now we prove the lower bound. Let $d$ be the integer such that
  $d-3<\eps\leq d-2$. Thus $d\geq3$. For all $n>3d$ such that $nd$ is
  even, \citet{Chandran03} constructed a graph $G$ with $n$ vertices,
  average degree $d\geq2+\eps$, and girth at least $(\log_d n)-1$. Now
  $G$ has Euler genus $g\leq\frac{dn}{2}-n+1\leq dn-2$. Thus $G$ has
  girth at least $(\log_d\frac{g+2}{d})-1$. Since $d<3+\eps$, the
  girth of $G$ is at least $h(\eps)\cdot \log(g+2)$ for some function
  $h$.
\end{proof}

We now extend Lemma~\ref{AverageDegreeFourPlusSees} for sufficiently
large embedded graphs.

\begin{lemma}
  \label{lem:Surface}
  Let $\eps\in(0,2)$.  Let $G$ be a graph with minimum degree 3 and
  average degree at least $4+\eps$. Assume that $G$ is embedded in a
  surface with Euler genus $g$, such that
  $|G|\geq(\frac{24}{\eps}+6)g$.  Then $G$ has a vertex $v$ that sees
  at most $2+\ceilFrac{12}{\eps}$ other vertices.
\end{lemma}

\begin{proof}
  For each vertex $v$ of $G$, associate a charge of
$$(8+2\eps)-\Brac{8+3\eps}\deg(v) + \Brac{24+6\eps}\frac{g}{|G|}
+\Brac{24+6\eps}\sum_{f\in F(G,v)}\frac{1}{|f|}\enspace.$$ Thus the
total charge is
\begin{align*}
  \;&(8+2\eps)|G|-\Brac{16+6\eps}\NumEdges{G} + \Brac{24+6\eps}g
  +\Brac{24+6\eps}|F(G)|\\
  =\;&(8+2\eps)|G|-\Brac{16+6\eps}\NumEdges{G} + \Brac{24+6\eps}g
  +\Brac{24+6\eps}(\NumEdges{G}-|G|-g+2)\\
  =\;&4(2\NumEdges{G} -(4+\eps)|G|) +
  2\Brac{24+6\eps}\\
  \geq \;&2\Brac{24+6\eps} \enspace.
\end{align*}
Thus some vertex $v$ has positive charge. That is,
\begin{align*}
  (8+2\eps)-\Brac{8+3\eps}\deg(v) + \Brac{24+6\eps}\frac{g}{|G|}
  +\Brac{24+6\eps}\sum_{f\in F(G,v)}\frac{1}{|f|} >0 \enspace.
\end{align*}
Since $\frac{(24+6\eps)g}{|G|} \leq \eps$,
\begin{align*}
  (24+6\eps) \sum_{f\in F(G,v)}\frac{ 1}{|f|} > ( 8+3\eps) (\deg(v)-1)
  \enspace.
\end{align*}
That is,
\begin{align*}
  \sum_{f\in F(G,v)}\frac{ 1}{|f|} >
  \Brac{\frac{1}{3}+\frac{1}{\alpha}} (\deg(v) -1) \enspace,
\end{align*}
where $\alpha:=6+\frac{24}{\eps}$.  We have proved that $\deg(v)$ and
the lengths of the facial walks incident to $v$ satisfy
Lemma~\ref{lem:NewTechnical} in Appendix~\ref{app}.  Thus $$\sum_{f\in
  F(G,v)}(|f|-2) \leq \left\lceil \frac{\alpha}{2} \right\rceil -1 =
2+\Ceil{\frac{12}{\eps}} \enspace.$$ The result follows from
Lemma~\ref{lem:Sees}.
\end{proof}

We now prove that the assumption that $n\in\Omega(\frac{g}{\eps})$ in
Lemma~\ref{lem:Surface} is needed.  Assume we are given $\eps\in(0,1]$
such that $k:=\frac{1}{\eps}-1$ is an integer. Hence $k\geq0$.
Consider the graph $G$ shown in Figure~\ref{fig:Surface}(b) with $2g$
vertices of degree $5$ and $2gk$ vertices of degree $4$. Thus
$|G|=2g(k+1)$ and
$2\NumEdges{G}=10g+8gk=2g(5+4k)=\frac{|G|}{k+1}(4k+5)=(4+\frac{1}{k+1})|G|=(4+\eps)|G|$.
Thus $G$ has average degree $4+\eps$.  Observe that every vertex lies
on a single face. Thus each vertex sees $|G|-1=\frac{2g}{\eps}-1$
other vertices.

A \emph{$k$-noose} in an embedded graph $G$ is a noncontractible
simple closed curve in the surface that intersects $G$ in exactly $k$
vertices. The \emph{facewidth} of $G$ is the minimum integer $k$ such
that $G$ contains a $k$-noose.

\begin{theorem}
  \label{th:SurfaceMinor}
  Let $\eps>0$.  Let $G$ be a 3-connected graph with average degree at
  least $4+\eps$, such that $G$ has an embedding in a surface with
  Euler genus $g$ and with facewidth at least $3$.  Then $G$ contains
  a $K_4$-model with at most $q(\eps)\cdot\log(g+2)$ vertices, for
  some function $q$. Moreover, for fixed $\eps$, this bound is within
  a constant factor of being optimal.
\end{theorem}

\begin{proof}
  If $|G|\leq(\frac{24}{\eps}+6)g$ then the result follows from
  Theorem~\ref{th:K4}. Otherwise, by Lemma~\ref{lem:Surface} some
  vertex $v$ sees at most $2+\ceilFrac{12}{\eps}$ other vertices.  The
  graph $G-v$ is 2-connected and has facewidth at least $2$.  Thus
  every face of $G-v$ is a simple cycle \citep[Proposition
  5.5.11]{MoharThom}.  In particular, the face of $G-v$ that contains
  $v$ in its interior is bounded by a simple cycle $C$. The vertices
  in $C$ are precisely the vertices that $v$ sees in $G$. Thus
  $G[C\cup\{v\}]$ is a subdivided wheel with $\deg(v)\geq3$
  spokes. Hence $G$ contains a $K_4$-model with at most
  $2+\ceilFrac{12}{\eps}$ vertices, which is at most
  $q(\eps)\cdot\log(g+2)$ for an appropriate function $q$.

  Now we prove the lower bound. Let $d$ be the integer such that
  $d-5<\eps\leq d-4$. Thus $d\geq5$. For every integer $n>3d$ such
  that $nd$ is even, \citet{Chandran03} constructed a graph $G$ with
  $n$ vertices, average degree $d\geq4+\eps$, and girth greater than
  $(\log_d n)-1$. Thus $G$ has Euler genus $g\leq\frac{dn}{2}\leq
  dn-2$.  Since every $K_4$-model contains a cycle, every $K_4$-model
  in $G$ has at least $(\log_d n)-1$ vertices. Since $n\geq
  \frac{g+2}{d}$ and $d<5+\eps$, every $K_4$-model in $G$ has at least
  $q(\eps)\cdot\log(g+2)$ vertices, for some function $q$.
\end{proof}

For a class of graphs, an edge is `light' if both its endpoints have
bounded degree.  For example, \citet{Wernicke-1904} proved that every
planar graph with minimum degree $5$ has an edge $vw$ such that
$\deg(v)+\deg(w)\leq 11$; see
\citep{Borodin-JRAM89,Kotzig55,JenMad96,JenVoss05} for extensions.
For a class of embedded graphs, we say an edge is `blind' if both its
endpoints see a bounded number of vertices.  In a triangulation, a
vertex only sees its neighbours, in which case the notions of `light'
and `blind' are equivalent.  But for non-triangulations, a `blind
edge' theorem is qualitatively stronger than a `light edge' theorem.
Hence the following result is a qualitative generalisation of the
above theorem of \citet{Wernicke-1904} (and of
Lemma~\ref{MinimumDegreeFiveSees}), and is thus of independent
interest. No such result is possible for minimum degree 4 since every
edge in $C_{2n}^2$ sees at least $n$ vertices.

\begin{proposition}
  \label{BlindEdgeSurface}
  Let $G$ be a graph with minimum degree $5$ embedded in a surface
  with Euler genus $g$, such that $|G|\geq240 g$. Then $G$ has an edge
  $vw$ such that $v$ and $w$ each see at most $12$ vertices. Moreover,
  for plane graphs (that is, $g=0$), $v$ and $w$ each see at most $11$
  vertices.
\end{proposition}

\begin{proof}
  Consider each vertex $x$. Let $\ell_x$ be the maximum length of a
  facial walk containing $x$. Let $t_x$ be the number of triangular
  faces incident to $x$, unless every face incident to $x$ is
  triangular, in which case let $t_x:=\deg(x)-1$.  Say $x$ is
  \emph{good} if $x$ sees at most $12$ vertices, otherwise $x$ is
  \emph{bad}.  Let
$$c_x:= 240-120\deg(x) + 240\frac{g}{|G|}+240 \sum_{f\in F(G,x)}\frac{1}{|f|}$$ be the
charge at $x$.  ($c_x$ is 240 times the \emph{combinatorial curvature}
at $x$.)\ By Euler's formula, the total charge is
$$240(|G| -
\NumEdges{G} + g + |F(G)|)=480\enspace.$$ Observe that (since
$\ell_x\geq3$ and $t_x\leq\deg(x)-1$ and $\deg(x)\geq5$)
\begin{align}
  c_x  \leq\;& 240-120\deg(x) + 240\frac{g}{|G|} + 240 \Brac{\frac{1}{\ell_x}+\frac{t_x}{3}+\frac{\deg(x)-t_x-1}{4}}\nonumber\\
  \leq\;& 181 -60\deg(x) + \frac{240}{\ell_x}+20t_x \label{eq:BlindEdgeSurface1}\\
  \leq\;& 241 -40\deg(x) \leq 41 \label{eq:BlindEdgeSurface2}\enspace.
\end{align}
For each good vertex $x$, equally distribute the charge on $x$ to its
neighbours. (Bad vertices keep their charge.)\ Let $c'_x$ be the new
charge on each vertex $x$.  Since the total charge is positive,
$c'_v>0$ for some vertex $v$.  If $v$ is good, then all the charge at
$v$ was received from its neighbours during the charge distribution
phase, implying some neighbour $w$ of $v$ is good, and we are
done. Now assume that $v$ is bad. Let $D_v$ be the set of good
neighbours of $v$. By \eqref{eq:BlindEdgeSurface1} and
\eqref{eq:BlindEdgeSurface2}, and since $\deg(w)\geq5$,
\begin{align}
  \label{eq:basicGenus}
  0<c'_v = c_v+\sum_{w\in D_v}\frac{c_w}{\deg(w)} \leq 181 -60\deg(v)
  + \frac{240}{\ell_v}+ 20t_v + \frac{41}{5}|D_v| \enspace.
\end{align}

We may assume that no two good neighbours of $v$ are on a common
triangular face.

\begin{claim}
  \label{cl:Claim}
  $|D_v|\leq\deg(v)-\frac{t_v}{2}$. Moreover, if
  $|D_v|=\deg(v)-\frac{t_v}{2}$ then some face incident to $v$ is
  non-triangular, and for every bad neighbour $w$ of $v$, the edge
  $vw$ is incident to two triangular faces.
\end{claim}

\begin{proof}
  First assume that every face incident to $v$ is triangular. Thus no
  two consecutive neighbours of $v$ are good. Hence
  $|D_v|\leq\frac{\deg(v)}{2}<\frac{\deg(v)+1}{2}=\deg(v)-\frac{t_v}{2}$,
  as claimed. This also proves that if $|D_v|=\deg(v)-\frac{t_v}{2}$
  then some face incident to $v$ is non-triangular.

  We prove the case in which some face incident to $v$ is
  non-triangular by a simple charging scheme. If $w$ is a good
  neighbour of $v$, then charge $vw$ by $1$. Charge each triangular
  face incident to $v$ by \half. Thus the total charge is
  $|D_v|+\frac{t_v}{2}$.  If $uvw$ is a triangular face incident to
  $v$, then at least one of $u$ and $w$, say $w$, is bad; send the
  charge of \half\ at $uvw$ to $vw$. Each good edge incident to $v$
  gets a charge of $1$, and each bad edge incident to $v$ gets a
  charge of at most \half\ from each of its two incident faces. Thus
  each edge incident to $v$ gets a charge of at most 1. Thus the total
  charge, $|D_v|+\frac{t_v}{2}$, is at most $\deg(v)$, as claimed.

  Finally, assume that $|D_v|=\deg(v)-\frac{t_v}{2}$. Then for every
  bad neighbour $w$ of $v$, the edge $vw$ gets a charge of exactly 1,
  implying $vw$ is incident to two triangular faces.
\end{proof}

Claim~\ref{cl:Claim} and \eqref{eq:basicGenus} imply
\begin{align*}
  \nonumber 0&< 181 -60\deg(v) + \frac{240}{\ell_v}+ 20t_v
  +\frac{41}{5}\deg(v)  - \frac{41t_v}{10} \\
  &= 181 -\frac{259}{5}\deg(v) + \frac{240}{\ell_v}+
  \frac{159}{10}t_v\enspace.
\end{align*}
Since $t_v\leq\deg(v)-1$ and $\deg(v)\geq5$,
\begin{align*}
  0 < \frac{1651}{10} -\frac{359}{10}\deg(v) + \frac{240}{\ell_v} \leq
  -\frac{144}{10} + \frac{240}{\ell_v} \enspace.
\end{align*}
implying $\ell_v\in\{3,4,\dots,16\}$. Since $\ell_v\geq3$,
\begin{align*}
  0 < \frac{2451}{10} -\frac{359}{10}\deg(v) \enspace,
\end{align*}
implying $\deg(v)\in\{5,6\}$ and $t_v\in\{0,1,\dots,\deg(v)-1\}$.

We have proved that finitely many values satisfy
\eqref{eq:basicGenus}.  We now strengthen this inequality in the case
that $|D_v|=\deg(v)-\frac{t_v}{2}$.

Let $f$ be a face of length $\ell_v$ incident to $v$. Let $x$ and $y$
be two distinct neighbours of $v$ on $f$. Suppose on the contrary that
$x$ is bad. By Claim~\ref{cl:Claim}, $vx$ is incident to two
triangular faces, one of which is $vxy$. Thus $\ell_v=3$, and every
face incident to $v$ is a triangle, which contradicts the Claim. Hence
$x$ is good. Similarly $y$ is good.

Thus $\ell_x\geq\ell_v$. By \eqref{eq:BlindEdgeSurface1},
\begin{align*}
  c_x \leq 181-60\deg(x) + \frac{240}{\ell_v}+20t_x \leq 161-40\deg(x)
  + \frac{240}{\ell_v} \leq \frac{240}{\ell_v} -39 \enspace.
\end{align*}
Similarly, $c_y\leq \frac{240}{\ell_v} -39$. Hence (assuming
$|D_v|=\deg(v)-\frac{t_v}{2}$),
\begin{align}
  \nonumber 0<c'_v &\leq 181-60\deg(v) + \frac{240}{\ell_v}+20t_v
  +\frac{c_x}{\deg(x)} +\frac{c_y}{\deg(y)}
  +\sum_{w\in D_v\setminus\{x,y\}}\frac{c_w}{\deg(w)}\\
  \nonumber &\leq 181-60\deg(v) + \frac{240}{\ell_v}+20t_v
  +\frac{\frac{240}{\ell_v}-39}{\deg(x)}
  +\frac{\frac{240}{\ell_v}-39}{\deg(y)}
  +\sum_{w\in D_v\setminus\{x,y\}}\frac{41}{\deg(w)}\\
  &\leq 181-60\deg(v) + \frac{240}{\ell_v}+20t_v
  +2\Brac{\frac{48}{\ell_v}-\frac{39}{5}}
  +\frac{41}{5}(|D_v|-2)\label{eq:MoreoverGenus} \enspace.
\end{align}

Checking all values of $\deg(v)$, $t_v$ and $\ell_v$ that satisfy
\eqref{eq:basicGenus} and \eqref{eq:MoreoverGenus} proves that
$$t_v+(\deg(v)-t_v)(\ell_v-2)\leq 12$$ (which is tight for $\deg(v)=5$
and $t_v=4$ and $\ell_v=10$ and $|D_v|=2$).  Thus $$\sum_{f\in
  F(G,v)}(|f|-2)\leq t_v(3-2)+(\deg(v)-t_v)(\ell_v-2)\leq
12\enspace.$$ By Lemma~\ref{lem:Sees}, $v$ sees at most $12$
vertices. Therefore $v$ is good, which is a contradiction.

In the case of planar graphs, we define a vertex to be \emph{good} if
it sees at most 11 other vertices. Since $g=0$,
\eqref{eq:BlindEdgeSurface1} and \eqref{eq:BlindEdgeSurface2} can be
improved to
\begin{align}
  c_x \leq 180 -60\deg(x) + \frac{240}{\ell_x}+20t_x \leq 240
  -40\deg(x) \leq 40\enspace.
\end{align}
Subsequently, \eqref{eq:basicGenus} is improved to
\begin{align}
  \label{eq:basicPlanar}
  0<c'_v = 180 -60\deg(v) + \frac{240}{\ell_v}+ 20t_v + 8|D_v|
  \enspace,
\end{align}
and \eqref{eq:MoreoverGenus} is improved to
\begin{align}
  \label{eq:MoreoverPlanar}
  0<c'_v \leq 180-60\deg(v) + \frac{240}{\ell_v}+20t_v
  +2\Brac{\frac{48}{\ell_v}-8} +8(|D_v|-2)\enspace.
\end{align}
Checking all values of $\deg(v)$, $t_v$ and $\ell_v$ that satisfy
\eqref{eq:basicPlanar} and \eqref{eq:MoreoverPlanar} proves that
$t_v+(\deg(v)-t_v)(\ell_v-2)\leq 11$. As in the main proof, it follows
that $v$ is good.
\end{proof}

\begin{figure}[!ht]
  \begin{center}
    \includegraphics{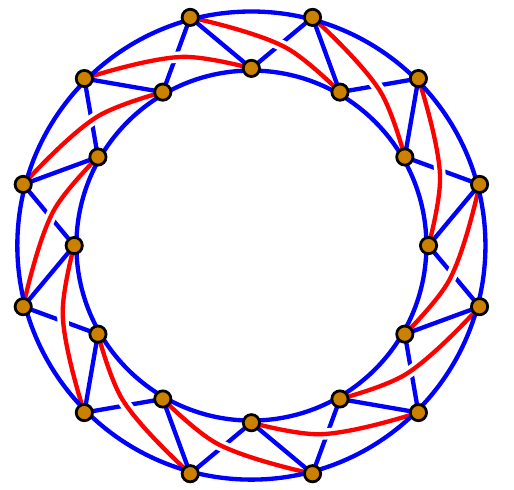}
  \end{center}
  \caption{\label{fig:CycleCubed}$C^2_{24}$ plus a perfect matching,
    embedded on $\mathbb{S}_{12}$.}
\end{figure}

We now prove that the assumption that $|G|\in\Omega(g)$ in
Proposition~\ref{BlindEdgeSurface} is necessary. Let $G$ be the graph
obtained from $C_{2n}^2$ by adding a perfect matching, as shown
embedded in $\mathbb{S}_n$ in Figure~\ref{fig:CycleCubed} (where there
is one handle for each pair of crossing edges). This graph is
5-regular, but each vertex is on a facial walk of length $n$. Thus no
vertex sees a bounded number of vertices.

\section{Open Problems}

The first open problem that arises from this work is to determine the
best possible function $f$ in Theorem~\ref{th:Main}. In particular,
does average degree at least some polynomial in $t$ force a small
$K_t$-model? Even stronger, is there a function $h$, such that every
graph $G$ with average degree at least $f(t) + \eps$ contains a
$K_t$-model with $h(t,\eps)\cdot\log |G|$ vertices, where $f(t)$ is
the minimum number such that every graph with average degree at least
$f(t)$ contains a $K_t$-minor? We have answered this question in the
affirmative for $t\leq 4$. The case $t=5$ is open. It follows from
Wagner's characterisation of graphs with no $K_5$-minor that average
degree at least $6$ forces a $K_5$-minor \citep{Wagner37}.
Theorem~\ref{th:Kt} proves that average degree at least $16+\eps$
forces a $K_5$-model with at most $h(\eps) \cdot \log n$ vertices.  We
conjecture the following improvement:

\begin{conjecture}
  There is a function $h$ such that for all $\eps>0$, every graph $G$
  with average degree at least $6+\eps$ contains a $K_5$-model with at
  most $h(\eps)\cdot\log|G|$ vertices.
\end{conjecture}

\begin{figure}[!ht]
  \begin{center}
    \includegraphics{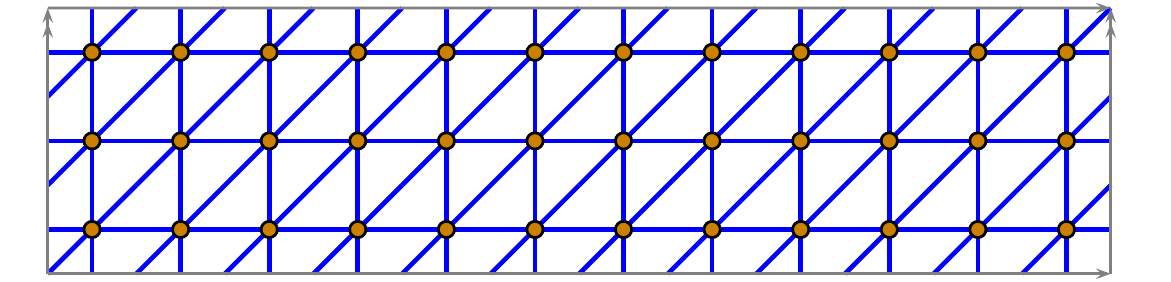}
  \end{center}
  \caption{\label{fig:Torus} 6-regular $12\times 3$ triangulated
    toroidal grid}
\end{figure}

This degree bound would be best possible: Let $G_n$ be the 6-regular
$n\times 3$ triangulated toroidal grid, as illustrated in
Figure~\ref{fig:Torus}. Every $K_5$-model in $G_n$ intersects every
column (otherwise $K_5$ is planar). Thus every $K_5$-model in $G_n$
has at least $n$ vertices.

Note that, while in this paper we have only studied small
$K_t$-models, the same questions apply for small $H$-models, for
arbitrary graphs $H$. This question was studied for $H=K_4-e$ in
\citep{FJP}. See \citep{Thomason06,MT-Comb05,Myers-DM03,KP08,KO05} for
results about forcing $H$-minors.

\subsection*{Acknowledgments}

Thanks to Michele Conforti for suggesting to study the relationship
between average degree and small models. Thanks to Paul Seymour for
suggesting the example following Lemma~\ref{le:NiceK3}. Thanks to
Alexandr Kostochka for pointing out reference \citep{KP88}.


\def\cprime{$'$}
\def\soft#1{\leavevmode\setbox0=\hbox{h}\dimen7=\ht0\advance \dimen7
  by-1ex\relax\if t#1\relax\rlap{\raise.6\dimen7
    \hbox{\kern.3ex\char'47}}#1\relax\else\if T#1\relax
  \rlap{\raise.5\dimen7\hbox{\kern1.3ex\char'47}}#1\relax \else\if
  d#1\relax\rlap{\raise.5\dimen7\hbox{\kern.9ex
      \char'47}}#1\relax\else\if D#1\relax\rlap{\raise.5\dimen7
    \hbox{\kern1.4ex\char'47}}#1\relax\else\if l#1\relax
  \rlap{\raise.5\dimen7\hbox{\kern.4ex\char'47}}#1\relax \else\if
  L#1\relax\rlap{\raise.5\dimen7\hbox{\kern.7ex
      \char'47}}#1\relax\else\message{accent \string\soft \space #1
    not defined!}#1\relax\fi\fi\fi\fi\fi\fi}

\appendix
\section{Some Technicalities}
\label{app}

\begin{lemma}
  \label{lem:technical}
  Let $\alpha>0$.  Let $d,f_1,\dots,f_d$ be integers, each at least
  $3$, such that
$$\sum_{i=1}^d\frac{1}{f_i}>\Brac{\frac{1}{3}+\frac{1}{\alpha}} d -\frac{1}{3}\enspace.$$
Then
$$\sum_{i=1}^d (f_i-2) \leq \left\lceil \frac{\alpha}{3} \right\rceil -1\enspace.$$
\end{lemma}

\begin{proof} 
We may assume that $f_1, \ldots, f_d$ firstly maximise $\sum_i (f_i
- 2)$, and secondly maximise $\sum_i \frac{1}{f_i}$. We claim that
$f_i = 3$ for all $i \in \{1,\ldots,d\}$ except perhaps one. Suppose
on the contrary that $f_j \geq f_k \geq 4$ for distinct $j, k \in
\{1,\ldots,d\}$. Let $f'_i := f_i$ for $i \in \{1,\ldots,d\} \setminus\{j,k\}$, $f'_j := f_j + 1$, and $f'_k := f_k - 1$. Then
$$
\sum_{i=1}^d f'_i = \sum_{i=1}^d f_i\,\quad \text{but}\quad
\sum_{i=1}^d \frac{1}{f'_i} > \sum_{i=1}^d \frac{1}{f_i}\enspace,
$$
implying $f_1,\dots,f_d$ do not maximise $\sum_j\frac{1}{f_j}$.  
Thus the claim holds and we may assume $f_i = 3$ for $i \in \{1,\ldots,d-1\}$. 
Hence
$$
\frac{d-1}{3} + \frac{1}{f_d} > \left( \frac{1}{3} + \frac{1}{\alpha} \right) d - \frac{1}{3}\enspace.
$$
Thus $\frac{1}{f_d}>\frac{d}{\alpha}$, implying $f_d \leq \ceil{ \frac{\alpha}{d} } - 1$. 
Since $\frac{\alpha}{d}>f_d\geq 3$ and since $d\geq3$,
$$
\frac{\alpha}{3}
= \frac{\alpha}{d} \Brac{ \frac{d}{3} - 1 } + \frac{\alpha}{d} 
\geq 3 \Brac{ \frac{d}{3} -  1 }  + \frac{\alpha}{d} 
= d-3 + \frac{\alpha}{d} \enspace.
$$
Hence
$$
\Ceil{\frac{\alpha}{3}} \geq\Ceil{ d - 3 + \frac{\alpha}{d}} = d - 3 + \Ceil{ \frac{\alpha}{d}}\enspace.
$$
Therefore
$$\sum_{i=1}^d(f_i-2) 
\leq  (d-1)(3-2)+ \Ceil{\frac{\alpha}{d}}-3
= d-3+ \Ceil{\frac{\alpha}{d}}-1 
\leq \CeilFrac{\alpha}{3}-1\enspace.$$
This completes the proof.
\end{proof}

\begin{lemma}
  \label{lem:NewTechnical}
  Let $\alpha>0$. Let $d,f_1,\dots,f_d$ be integers, each at least
  $3$, such that
$$\sum_{i=1}^d\frac{1}{f_i}>\Brac{\frac{1}{3}+\frac{1}{\alpha}} ( d -1)\enspace.$$
Then
$$\sum_{i=1}^d (f_i-2) \leq \CeilFrac{\alpha}{2} -1 \enspace.$$
\end{lemma}

\begin{proof}   
  As in the proof of Lemma~\ref{lem:technical}, we may assume that
  $f_j= 3$ for all $j\in\{3,\dots,d-1\}$.  Hence
$$\frac{d-1}{3}+\frac{1}{f_d}>\Brac{\frac{1}{3}+\frac{1}{\alpha}} ( d -1)\enspace.$$
Thus $\frac{1}{f_d}>\frac{d-1}{\alpha}$, implying $f_d \leq \ceil{
  \frac{\alpha}{d-1} } - 1$.  Since $\frac{\alpha}{d-1}>f_d\geq 3$ and
since $d\geq3$,
$$
\frac{\alpha}{2} \geq \frac{\alpha d}{3(d-1)} =
\Brac{\frac{\alpha}{d-1}} \Brac{ \frac{d}{3} - 1 } +
\frac{\alpha}{d-1} \geq 3 \Brac{ \frac{d}{3} - 1} + \frac{\alpha}{d-1}
= d-3 + \frac{\alpha}{d-1}\enspace.
$$
Hence
$$
\Ceil{\frac{\alpha}{2}} \geq\Ceil{ d - 3 + \frac{\alpha}{d-1}} = d - 3
+ \Ceil{ \frac{\alpha}{d-1}}\enspace.
$$
Therefore
$$\sum_{i=1}^d(f_i-2) 
\leq (d-1)(3-2)+ \Ceil{\frac{\alpha}{d-1}}-3 = d-3+
\Ceil{\frac{\alpha}{d-1}}-1 \leq \Ceil{ \frac{\alpha}{2}}-1
\enspace.$$ This completes the proof.
\end{proof}

\end{document}